\documentclass[12pt]{article}
\usepackage[dvips]{graphicx}
\usepackage{amsmath,amssymb}
\usepackage{bm}

\title{Brief note on Riemann hypothesis}
\author{}

\begin{document}
{\pagestyle{empty}
\rightline{June 2009}
\rightline{~~~~}
\vskip 1cm
\centerline{\large \bf A Brief note on the Riemann hypothesis}
\vskip 1cm
\centerline
{{Minoru Fujimoto
 } and
 {Kunihiko Uehara
 }
}
\vskip 1cm
\centerline{\it ${}^1$Seika Science Research Laboratory,
Seika-cho, Kyoto 619-0237, Japan}
\centerline{\it ${}^2$Department of Physics, Tezukayama University,
Nara 631-8501, Japan}
\vskip 2cm

\centerline{\bf Abstract}
\vskip 0.2in
  We deal with the Euler's alternating series 
of the Riemann zeta function 
to define a regularized ratio appeared in the functional equation
even in the critical strip and 
show some evidence to indicate the hypothesis in this note. 

\vskip 0.4cm\noindent

\noindent
PACS number(s): 02.30.-f, 02.30.Gp, 05.40.-a

\hfil
}
\setcounter{equation}{0}
\addtocounter{section}{0}
\hspace{\parindent}

  Regularizations by way of the zeta function have been successful with 
some physical applications so far, but the Riemann hypothesis associated with the Riemann zeta 
function itself has been remained to be proved. 
Recently we have proposed a regularization technique \cite{Fujimoto1,Fujimoto2} and apply this 
regularization to the Euler product of zeta functions. 
Here we utilize the Euler's alternating summation, 
which is finite even in the critical strip and 
seems to be essential to clarify the Riemann hypotheses.

  The definition of the Riemann zeta function is 
\begin{equation}
  \zeta(z)=\lim_{n\to\infty}\zeta_n(z),\ 
  \zeta_n(z)\equiv\sum_{k=1}^n\frac{1}{k^z}
\label{e101}
\end{equation}
for $\Re z>1$. 
In this note we adopt a hat notation like $\hat{\zeta}(z)$ for $\Re z>0$ such as 
\begin{equation}
  \hat\zeta(z)=
               \frac{1}{1-2^{1-z}}\sum_{n=1}^\infty\frac{(-1)^{n-1}}{n^z},
\label{e102}
\end{equation}
which is well defined even in the critical strip $0<\Re z<1$ and 
a summation part is called as the Euler's alternating series. 
We often express a hat representation by a ``regularized" form because 
a hat representation is defined by a subtraction 
an infinite number from a divergent quantity.

  In this note, we deal with the Euler's alternating series 
of the Riemann zeta function as (\ref{e102}) 
to well-define even in the critical strip $0<\Re z<1$ 
and utilize the functional equation to indicate the hypothesis. 
Hereafter we are only interested in the region $\Re z\ge\frac{1}{2}$ 
for the Riemann zeta function, because the functional equation ensures 
the regularized nature of the zeta function for the other half plane 
$\Re z<\frac{1}{2}$.

\vskip 5mm
  There is a relation called the functional equation of the Riemann zeta function 
\begin{equation}
  \hat{\zeta}(z)=\hat{H}(z)\hat{\zeta}(1-z),
\label{e121}
\end{equation}
where $\hat{H}(z)$ is given by $\displaystyle 2\Gamma(1-z)(2\pi)^{z-1}\sin\frac{\pi z}{2}$, 
which we deal with for the infinite limit of $\hat{H}_n(z)$ defined by 
\begin{equation}
  \hat{H}_n(z)=\frac{\hat{\zeta}_n(z)}{\hat{\zeta}_n(1-z)},
\label{e122}
\end{equation}
where $\hat{\zeta}_n(z)$ is defined by
\begin{equation}
  \hat{\zeta}_n(z)\equiv\zeta_n(z)-\frac{n^{1-z}}{1-z}.
\label{e123}
\end{equation}

  When we think about the Euler's alternating series in (\ref{e102}) 
for the Riemann zeta function
\begin{equation}
  \hat{\zeta}(z)=\lim_{n\to\infty}\frac{1}{1-2^{1-z}}\xi_n(z),
\label{e201}
\end{equation}
where $\displaystyle\xi_n(z)\equiv\sum_{k=1}^n\frac{(-1)^{k-1}}{k^z}$, 
we can evaluate the function for $z$ even in the critical strip $0<\Re z<1$ 
as mentioned above. 
The relation between the zeta function and $\xi_n(z)$ is special 
because the relation form itself conserves before and after the regularization 
as
\begin{eqnarray}
  \xi_{2n}(z)&=&\zeta_{2n}(z)-2^{1-z}\zeta_n(z),
\label{e202}\\
  \xi_{2n}(z)&=&\hat{\zeta}_{2n}(z)-2^{1-z}\hat{\zeta}_n(z),
\label{e203}
\end{eqnarray}
where we used the relation (\ref{e123}).

  When we put $z=\rho$ and take the limit of $n\to\infty$ in (\ref{e203}), 
where $\rho$ is one of the non-trivial zeroes for the Riemann zeta function, 
we get
\begin{equation}
  \lim_{n\to\infty}\hat{\zeta}_{2n}(\rho)=2^{1-\rho}\lim_{n\to\infty}\hat{\zeta}_n(\rho)
\label{e204}
\end{equation}
and using the property that $1-\rho$ is also a zero as $\rho$ is, we also get
\begin{equation}
  \lim_{n\to\infty}\hat{\zeta}_{2n}(1-\rho)=2^{\rho}\lim_{n\to\infty}\hat{\zeta}_n(1-\rho).
\label{e205}
\end{equation}
  Combining (\ref{e204}) with (\ref{e205}), we get 
\begin{equation}
  \lim_{n\to\infty}\hat{H}_{2n}(\rho)=2^{1-2\rho}\lim_{n\to\infty}\hat{H}_n(\rho). 
\label{e206}
\end{equation}
Repeatedly applying (\ref{e206}) $m$-times to itself, we get
\begin{equation}
  \lim_{n\to\infty}\hat{H}_{2^m n}(\rho)=2^{m(1-2\rho)}\lim_{n\to\infty}\hat{H}_n(\rho).
\label{e208}
\end{equation}

When we think about the limit of $m\to\infty$ in (\ref{e208}), 
the left hand side will coincide with $\hat{H}(\rho)$ which is finite 
and $\displaystyle\lim_{n\to\infty}\hat{H}_n(\rho)$ is also finite in right hand side. 
After all we can conclude that the term $\displaystyle\lim_{m\to\infty}2^{m(1-2\rho)}$ is finite
which means that a real part of the zero $\Re\rho$ is identical to one half. 


  Another way to reach the conclusion is taking an absolute value 
in (\ref{e206}) as 
\begin{equation}
  |2^{1-2\rho}|=\lim_{n\to\infty}\left|\frac{\hat{H}_{2n}(\rho)}{\hat{H}_n(\rho)}\right|
               =\frac{|\hat{H}(\rho)|}{|\hat{H}(\rho)|}=1,
\label{e209}
\end{equation}
which again claims that a real part of the zero $\Re\rho$ is equal to 
$\displaystyle\frac{1}{2}$. 

\vskip 5mm
  Here we refer to an error estimation in (\ref{e204}) and (\ref{e205}) as 
\begin{eqnarray}
  \hat{\zeta}_{2n}(\rho)  &=&2^{1-\rho}\hat{\zeta}_n(\rho)+o(n),\\
  \hat{\zeta}_{2n}(1-\rho)&=&2^{\rho}\hat{\zeta}_n(1-\rho)+o(n),
\end{eqnarray}
which enable us to evaluate (\ref{e206}) or (\ref{e208}),
are found by making use of the relation derived from Hardy and Littlewood
\begin{equation}
  \hat{\zeta}(z)=\hat{\zeta}_n(z)+O(n^{-\Re z})
\end{equation}
for $|\Im z|\le 2\pi n/C$, where $C$ is constant greater than one.

  Finally we have to mention that the regularized form for the Riemann zeta function above 
in the region $\frac{1}{2}\le\Re z<1$ coincides with one by the analytic continuation.

\vskip 5mm

\end{document}